\documentclass{article}
\usepackage{amssymb}
\usepackage{amsmath}
\usepackage{amsthm}
\usepackage{amsbsy}
\usepackage{psfrag}
\usepackage{graphics}
\usepackage{graphicx}
\usepackage{subfigure}
\usepackage{epsfig}
\usepackage{epstopdf}
\usepackage{subfig}

\theoremstyle{definition}

\theoremstyle{remark}

\numberwithin{equation}{section}

\newcommand{\x}{\textbf{x}}
\newcommand{\y}{\textbf{y}}
\usepackage{graphicx}
\usepackage{psfrag}
\usepackage{caption}
\usepackage{float}

\begin{document}

\title{\bf Double Negative Dispersion Relations from Coated Plasmonic Rods}

\author{Yue Chen \\ {\normalsize Department of Mathematics} \\
{\normalsize Louisiana State University}\\
{\normalsize Baton Rouge, LA 70803, USA.}\\ {\normalsize email: chenyue\symbol{'100}math.lsu.edu}\\[3pt]
\and Robert Lipton \\{\normalsize Department of Mathematics}\\{\normalsize  Louisiana State University}\\
{\normalsize Baton Rouge, LA 70803, USA.}\\ {\normalsize email: lipton\symbol{'100}math.lsu.edu}\\[3pt]}

\maketitle

\begin{abstract}

A metamaterial with frequency dependent double negative effective properties is constructed from a sub-wavelength periodic array of coated rods. Explicit power series are developed for the dispersion relation and associated Bloch wave solutions. The expansion parameter is the ratio of the length scale of the periodic lattice to the wavelength.  Direct numerical simulations for finite size period cells show that the leading order term in the power series for the dispersion relation is a good predictor of the dispersive behavior of the metamaterial. 

\bigskip


\end{abstract}

\section{Introduction}
\label{introduction}
\par
Metamaterials are artificial materials designed to have electromagnetic properties not generally found in nature. One contemporary area of research explores novel sub-wavelength constructions that deliver  metamaterials with both negative bulk dielectric constant and bulk  magnetic permeability across certain frequency intervals. These double negative materials are promising materials for the creation of negative index super lenses that overcome the small diffraction limit and have great potential in applications such as biomedical imaging, optical lithography and data storage. The  early work of Veselago \cite{Veselago} identified novel effects associated with hypothetical  materials for which both the  dielectric constant and magnetic permeability  are simultaneously negative. Such double negative media support electromagnetic wave propagation in which the phase velocity is antiparallel to the direction of energy flow, and other unusual electromagnetic effects such as the reversal of the Doppler effect and Cerenkov radiation. At the end of the last century  Pendry \cite{Pendry1998}  
demonstrated that unconventional properties can be derived from subwavelength configurations of different conventional materials and showed that a cubic lattice of metal wires exhibits behavior associated with negative bulk dielectric constant. Subsequently it was shown that  a periodic array of non-magnetic metallic split-ring resonators deliver negative effective magnetic permeability at microwave frequencies  \cite{PendryHolden}.  In more recent work  Smith et al. \cite{Smith} experimentally demonstrated that metamaterials made from arrays of metallic posts and split ring resonators generate an effective negative refractive index at microwave frequencies.  Building on this Shelby et al. \cite{ShelbySmith}  experimentally confirmed that a microwave beam would undergo negative refraction at the interface between such a metamaterial and air. Subsequent work has delivered several new designs using different configurations of metallic resonators for double negative behavior \cite{23,17,21,18,20,22}.

For higher frequencies in the infrared and optical range,  new strategies for generating materials with double negative bulk properties  rely on Mie resonances. One scheme employs coated rods made from a high dielectric core coated with a frequency dependent dielectric plasmonic or Drude type behavior at optical frequencies  \cite{11,Yannopappas,Yanno2}. A second scheme employs small rods or particles made from dielectric materials with large  permittivity, \cite{Plasmon, LPeng, VinkFelbacq}. Alternate strategies for generating negative bulk dielectric permeability at infrared and optical frequencies  use special configurations of plasmonic nanoparticles \cite{6A}, \cite{shevts}. The list of  metamaterial systems is rapidly growing  and comprehensive reviews of the subject can be found in \cite{Service} and \cite{Shalaev}.
\par
\begin{figure}[h!]
\begin{center}
\begin{psfrags} 
\psfrag{P}{$P$} 
\psfrag{R}{$R$}
\psfrag{H}{$H$}
\includegraphics[width=0.4\textwidth]{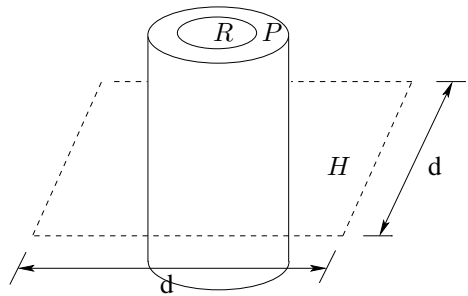}
\end{psfrags} 
\caption{Coated cylinder microgeometry: $R$ represents the high dielectric core, $P$ the plasmonic coating and $H$ denotes the connected host material. }
\label{coatedcylindergeo}
\end{center}
\end{figure}

In this article we construct metamaterials made from  subwavelength periodic arrangements of nonmagnetic  infinitely long coated cylinders immersed in a nonmagnetic host. 
The coated cylinders are parallel to the $x_3$ axis  and made from a  frequency independent high dielectric core and a frequency dependent dielectric plasmonic coating (Figure \ref{coatedcylindergeo}). 
We apply the mathematical analysis developed by the authors  in \cite{chenlipton} to  express the effective dielectric constant and magnetic permeability in terms of spectral representation formulas. These formulas are determined by the Dirichlet spectra of the core and the generalized electrostatic resonances associated with the region exterior to the core.  The formulas are used to calculate the frequency intervals where either double negative or double positive bulk properties appear. These intervals are governed by the poles and zeros of the effective magnetic permeability and effective dielectric permittivity tensors.  Explicit power series are developed for the dispersion relation and associated Bloch wave solutions. The frequency intervals over which the effective properties are either double negative or double positive imply the existence of convergent power series representations  for Bloch wave modes in the dynamic regime away from the homogenization limit. 
We apply the power series representation to calculate the average Poynting vector and show that in the homogenization limit the energy flow and phase velocity are in opposite directions over frequency intervals associated with double negative behavior.
This gives the requisite explicit and mathematically rigorous analysis beyond the homogenization limit  and provides evidence for wave propagation in the double negative regime for this class of metamaterial.

We apply the methods to a metamaterial made from a periodic array of circular coated cylinders making use of the method of Rayleigh \cite{Rayleigh} to numerically calculate the generalized  electrostatic resonances. These resonances together with the Dirichlet spectra of the core are used to identify explicit frequency intervals over which effective properties are double negative or double positive. Several branches of the  leading order dispersion relation are calculated using the spectral representation formulas for the effective magnetic permeablity and dielelectric constant. 
We compare these  with direct numerical simulations to find that the leading order dispersion relation is a good predictor of the dispersive behavior of the metamaterial. It is found that the leading order behavior trends with the direct numerical simulation  even when the length scale of the microstructure is only $20 \%$ smaller than the wavelength of the propagating wave. 
These results provide new methods necessary to identify frequency intervals characterized by negative index behavior and its influence on wave propagation beyond the homogenization limit.

Related work delivers formulas for frequency-dependent effective magnetic permeability together with conditions for generation of negative effective permeability \cite{bouchetteschwizer, felbacqbouchette, bouchettefelbacq, 4A, 9A, KohnShipman}.
For periodic arrays made from metal fibers a homogenization theory delivering negative effective dielectric constant \cite{bouchettebourel}  has been established. A novel method for creating metamaterials with prescribed effective dielectric permittivity and effective magnetic permeability at a fixed frequency
is developed in \cite{Milton2}.  New methodologies  for computing homogenized properties for metamaterials are presented in \cite{Alu}, \cite{SmithPendry}.

We conclude noting that the power series approach to sub-wavelength  analysis has been utilized and developed in \cite{FLS} for characterizing the  dynamic dispersion relations for Bloch waves inside plasmonic crystals. It has also been applied to assess the influence of effective negative permeability on the propagation of Bloch waves inside high contrast dielectrics \cite{FLS2},  the generation of negative permeability inside metallic - dielectric resonators \cite{shipman}, and for concentric coated cylinder assemblages generating a double negative media \cite{chenlipton2}.

\section{Power series representations}
\label{background}
\begin{figure}[h!]
\begin{center}
\begin{psfrags} 
\psfrag{x}{$x_1$} 
\psfrag{y}{$x_2$}
\psfrag{0.5}{$\frac{d}{2}$}
\psfrag{H}{$H$}
\psfrag{P}{$P$}
\psfrag{R}{$R$}
\psfrag{G}{$G$}
\psfrag{theta}{$\theta$}
\psfrag{a}{$a$}
\psfrag{b}{$b$}
\psfrag{r}{$r$}
\includegraphics[width=0.4\textwidth]{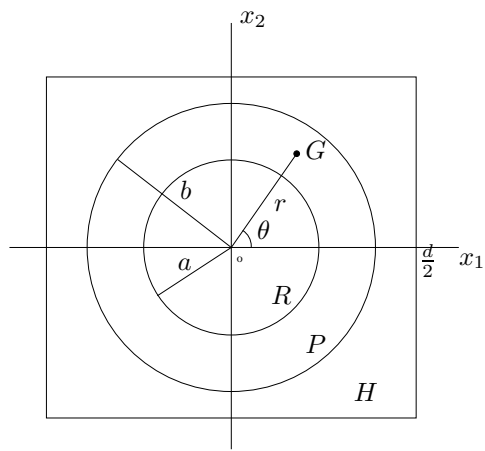}
\end{psfrags} 
\caption{The period cell}
\label{unitcell}
\end{center}
\end{figure}

We start with a metamaterial crystal characterized  by a period cell containing  a centered coated cylinder with plasmonic coating and high dielectric core. The core radius and the coating radius are denoted by $a$  and $b$ respectively (Figure \ref{unitcell}). 
The cylinder is parallel to the $x_3$ axis and is periodically arranged within a square lattice over the transverse $\mathbf{x}=(x_1,x_2)$ plane. The period of the lattice is denoted by $d$.  For  H-polarized Bloch-waves,  the magnetic field is aligned with the cylinders and the electric field lies in the transverse plane. The direction of propagation is described by the unit vector $\hat{\kappa}=(\kappa_1,\kappa_2)$ and $k=2\pi/\lambda$ is the wave number for a wave of length~$\lambda$ and the fields are of the form
\begin{eqnarray}
H_3=H_3(\mathbf{x})e^{i(k\hat{\kappa}\cdot\mathbf{x}-t\omega/c )},\,\,E_1=E_1(\mathbf{x})e^{i(k\hat{\kappa}\cdot\mathbf{x}-t\omega/c)},\,\,
E_2=E_2(\mathbf{x})e^{i(k\hat{\kappa}\cdot\mathbf{x}-t\omega/c)} \label{em3}
\end{eqnarray}
where $H_3(\mathbf{x})$, $E_1(\mathbf{x})$, and $E_2(\mathbf{x})$ are  $d$-periodic for $\mathbf{x}$ in $\mathbb{R}^2$. 
Here $c$ denotes the speed of light in free space. 
We denote the unit vector pointing along the $x_3$ direction by ${\bf e}_3$, and the periodic dielectric permittivity and magnetic permeability are denoted by $a_d$ and $\mu$ respectively. The electric field component $\mathbf{E}=(E_1,E_2)$ of the wave is determined by 
\begin{eqnarray}
 {\bf E}=-\frac{ic}{\omega a_d}{\bf e}_3\times \nabla H_3.
\label{relationofefandmf}
\end{eqnarray}

The materials are assumed non-magnetic hence the magnetic permeability $\mu$ is set to unity inside the coated cylinder and host. 
The oscillating dielectric permittivity for the crystal 
is a $d$  periodic function in the transverse plane and
is described by $a_d=a_d(\mathbf{x}/d)$ where $a_d(\mathbf{y})$ is 
the unit periodic dielectric function taking the values
\begin{equation}
a_d(\y)=
\begin{cases}
\epsilon_H &\text{ in the host material},\\
\epsilon_P(\omega) & \text{ in the frequency dependent ``plasmonic'' coating},\\
\epsilon_R=\epsilon_r/d^2 &\text{ in the high dielectric core} .
\end{cases}
\end{equation}
Here $\epsilon_r$ has dimensions of area and the frequency dependent permittivity $\epsilon_P$ of the plasmonic coating is given by \cite{11}, \cite{Yannopappas}  
\begin{eqnarray} 
\epsilon_P(\omega^2)=1-\frac{\omega_p^2}{\omega^2},
\label{singleosc}
\end{eqnarray} 
where $\omega$ is the frequency and $\omega_p$ is the plasma frequency \cite{bohren}. 
Setting $h^d(\mathbf{x})=H_3(\mathbf{x})e^{i(k\hat{\kappa}\cdot\mathbf{x})}$ the Maxwell equations take the form of the Helmholtz equation given by
\begin{equation}
-\nabla_{\mathbf{x}} \cdot \left(a_d^{-1}(\frac{\mathbf{x}}{d})\nabla_{\mathbf{x}} h^d(\mathbf{x})\right)=\frac{\omega^2}{c^2}h^d ~~~\text{ in  } \mathbb{R}^2.
\label{Helmholtzr2}
\end{equation}

We set $\x=d\y$ for $\y$ inside the unit period $Y=[-0.5,0.5]^2$, put  $\beta=dk\hat{\kappa}$ and write $u(\y)=H_3(d\y)$. The dependent variable is written $u^d(\y)=h^d(d\y)=u(\y)\exp^{i\beta\cdot\y}$,  and we recover the equivalent problem over the unit period cell given by
\begin{equation}
-\nabla_{\y} \cdot \left(a_d^{-1}(\y)\nabla_{\y} u^d\right)=\frac{d^2\omega^2}{c^2}u^d ~~~\text{ in  } Y.
\label{Helmholtzd}
\end{equation}
We start by introducing the power series in terms of dimensionless groups given by the ratio  $\rho=d/\sqrt{\epsilon_r}$, wave number $\tau=\sqrt{\epsilon_r}k$ and square frequency $\xi=\epsilon_r\frac{\omega^2}{c^2}$. The dimensionless parameter  measuring the departure away from  quasistatic or homogenization limits is given by the  ratio of period size to wavelength $\eta=dk=\rho\tau\geq 0$. 
For the problem considered here the $d=0$ limit is distinct from the quasi-static limit $k=0$. This is due to the explicit dependence of the dielectric constant on $d$ inside the core region of the rod.

For these parameters the dielectric permittivity takes the values  $\epsilon_P=1-\frac{\epsilon_r\omega_p^2/c^2}{\xi}$, $\epsilon_R=\frac{1}{\rho^2}$, $\epsilon_H=1$,  and  is denoted by $a_\rho(\y)$ for $\y$ in $Y$ and \eqref{Helmholtzd} is given by
\begin{equation}
-\nabla_{\y} \cdot \left(a_\rho^{-1}(\y)\nabla_{\y} u^d(\y)\right)=\rho^2\xi u^d(\y) ~~~\text{ in  } Y.
\label{Helmholtz}
\end{equation}
We introduce the space $H^1_{per}(Y)$ of trial and test functions that are square integrable with periodic boundary conditions on $Y$ and square integrable derivatives.  
The equivalent variational form of \eqref{Helmholtz} is given by
\begin{equation}
\int_Y a_\rho^{-1}\nabla u^d\cdot \nabla \bar{\tilde{v}}= \int_Y \frac{\rho^2\xi}{c^2}u^d\bar{\tilde{v}}
\label{variational}
\end{equation}
for any $\tilde{v}=v(\y)e^{i\hat{\kappa}\cdot \tau\rho \y}$ , where $v \in H ^1_{per}(Y)$ .

The unit period cell for the generic metamaterial system is represented in Figure \ref{unitcell}. In what follows $R$ represents the rod core cross section containing high dielectric material, $P$ the coating containing the plasmonic material and $H$ denotes the connected host material. Following \cite{chenlipton} we introduce the expansions for the Bloch wave $u^d$ eigenvalue $\xi$ pair
\begin{eqnarray}
 &&u^d=\sum_{m=0}^\infty \eta^m u_m e^{i\hat{\kappa}\cdot\tau\rho{\mathbf y}}\label{upower}
\\&&\xi=\sum_{m=0}^\infty\eta^m\xi_m \label{xipowerexpansion}
\end{eqnarray}
where $u_m$ belong to $H^1_{per}(Y)$. In view of the algebra it is convenient to write $u_m=i^m\underline{u}_0\psi_m$ where $\underline{u}_0$ is an arbitrary constant factor. Substitution of \eqref{upower} and \eqref{xipowerexpansion}   into \eqref{variational} and equating like powers of $\eta$ delivers an infinite coupled system of equations that can be solved iteratively.  We now describe this system of equations expressed in variational form.
Set $$z=\epsilon_P^{-1}(\xi_0)=\left(1-\frac{\epsilon_r\omega_p^2/c^2}{\xi_0}\right)^{-1}$$
and for $u$, $v$ belonging to $H^1_{per}(Y)$ we introduce the sesquilinear form
\begin{eqnarray}
B_z(u,v)=\int_{H}\nabla u\cdot\nabla\overline{v}\,d\y+\int_P\,z\,\nabla u\cdot\nabla\overline{v}\,d\y.
\label{bilinear}
\end{eqnarray}
Here $Y\setminus R=H\cup P$.
Substitution of the series 
into  (\ref{variational}) and equating like powers of $\eta$, produces the infinite set of coupled equations for $m=0,1,2\ldots$ given by
\begin{eqnarray}
&&\tau^2 B_z(\psi_m,v)\nonumber
\\&&+\xi_0^{-1}\epsilon_p^{-1}(\xi_0)\tau^2\int_{Y\setminus R}\big[\sum^{m-1}_{l=1}(-i)^l\xi_l\nabla \psi_{m-l}\cdot \nabla \overline{ v} +\hat{\kappa}\cdot\sum_{l=0}^{m-1}(-i)^l\xi_l(\psi_{m-1-l}\nabla \overline{v}-\nabla \psi_{m-1-l}\overline{v})\nonumber
\\&&-\sum_{l=0}^{m-2}(-i)^l\xi_l\psi_{m-2-l}\overline{v}\big] -\xi_0^{-1}\epsilon_p^{-1}(\xi_0)\tau^2\epsilon_r\frac{\omega_p^2}{c^2}\int_H\big[\hat{\kappa}\cdot(\psi_{m-1}\nabla \overline{v}-\nabla \psi_{m-1}\overline{v})-\psi_{m-2}\overline{v}\big]\nonumber
\\&&-\xi_0^{-1}\epsilon_p^{-1}(\xi_0)\int_R\big[\sum^{m-2}_{l=0}(-i)^l\xi_l\nabla \psi_{m-2-l}\cdot \nabla \overline{ v} +\hat{\kappa}\sum_{l=0}^{m-3}(-i)^l\xi_l(\psi_{m-3-l}\nabla \overline{v}-\nabla \psi_{m-3-l}\overline{v})\nonumber
\\&&-\sum_{l=0}^{m-4}(-i)^l\xi_l\psi_{m-4-l}\overline{v}\big] +\xi_0^{-1}\epsilon_p^{-1}(\xi_0)\int_R\epsilon_r\frac{\omega_p^2}{c^2}\big[\nabla \psi_{m-2}\cdot \nabla \overline{ v}+\hat{\kappa}(\psi_{m-3}\nabla \overline{v}-\nabla \psi_{m-3}\overline{v})+\psi_{m-4}\overline{v}\big]\nonumber
\\&&-\xi_0^{-1}\epsilon_p^{-1}(\xi_0)\int_Y\big[\sum_{l=0}^{m-2}\sum_{n=0}^l \xi_{m-2-l}\xi_n \psi_{l-n} i^{l-n-m}\overline{ v}+\epsilon_r\frac{\omega_p^2}{c^2}\sum_{l=0}^{m-2}(-i)^l\xi_l \psi_{m-2-l}\overline{ v}\big]\nonumber\\
&&=0,\hbox{ for all $v$ in $H^1_{per}(Y)$}.
\label{summation2}
\end{eqnarray}
Here the convention is $\psi_m=0$ for $m<0$.

Applying the theory developed in \cite{chenlipton}  the solvability of the infinite system for determining the unknown functions $\{\psi_m\}_{m=1}^\infty$ depends on the 
Dirichlet spectra of $R$ and the generalized electrostatic spectra associated with $Y\setminus R=P\cup H$. Here the Dirichlet spectra is  given by the eigenvalues $\mu_n>0$, $\mu_{n+1}\geq\mu_n$, $\mu_n\rightarrow \infty$, for $n\rightarrow\infty$ associated with the Dirichlet eigenfunctions of the Laplacian  on $R$. The generalized electrostatic spectra is characterized by all eigenvalues $\lambda$ and eigenfunctions $u$ of
\begin{equation}
\begin{cases}
\Delta u=0 ~~~\text{ in  } H ,\\
\Delta u=0 ~~~\text{ in  } P ,\\
\end{cases}
\end{equation}
with the boundary conditions
\begin{eqnarray}
 \begin{cases}
  u|^{-}=u|^{+} ~~~\text{ on } \partial P,\\
\partial_{r}u|_{r=a}=0 ~~~\text{ on } \partial R,\\
\lambda[\partial_r u]^{-}_{+}=-\frac{1}{2}(\partial_ru^-+\partial_ru^+) \text{ on } \partial P ,\\
u\text{ is $Y$-periodic } . \\
 \end{cases}
\label{boundarycondition}
\end{eqnarray}
The generalized electrostatic spectra is denumerable lies in the open interval $(-1/2,1/2)$ with zero being the only accumulation point, see \cite{chenlipton}. The eigenfunctions $\{\psi_{\lambda_n}\}_{n=0}^\infty$  associated with the electrostatic resonances $\{\lambda_n\}_{n=1}^\infty$ form a complete orthonormal set of functions in the space of  mean zero periodic functions belonging to $H^1_{per}(Y\setminus R)$ that are harmonic in $P$ and $H$, \cite{chenlipton}. Here orthonormality is with respect to the inner product $(u,v)=\int_{Y\setminus R}\nabla u\cdot\nabla\overline{v}\,dx$. The complete orthonormal systems of eigenfunctions associated with electrostatic resonances and Dirichlet eigenvalues are used to solve for $\psi_0$ and $\psi_1$ in $H\cup P$ and provide an explicit formula for $\xi_0$. We follow \cite{chenlipton} to find that $\psi_0=1$ in $Y\setminus R$ and 
\begin{eqnarray}
-\Delta\psi_0=\xi_0\psi_0, \hbox{ in $R$ } 
\label{Laplace}
\end{eqnarray}
with $\psi_0=1$ on the boundary of $R$.
From \eqref{summation2}   we find that  $\psi_1$ 
is the solution of 
\begin{eqnarray}
-\Delta \psi_1&=&0, \hbox{in $P$ and in $H$}
\label{deltah1}
\end{eqnarray}
and the corresponding transmission conditions for $\psi_1$  are given by
\begin{eqnarray}
n\cdot\left(\nabla \psi_1+ i\hat{\kappa} \right)_{\scriptscriptstyle_{|_H}}&=&n\cdot\epsilon_P^{-1}(\xi_0)\left(\nabla \psi_1 +i\hat{\kappa}\right)_{\scriptscriptstyle_{|_P}},\hbox{ H-P interface},
\label{HC1}\\
n\cdot\epsilon_P^{-1}(\xi_0)\left(\nabla \psi_1 +i\hat{\kappa} \right)_{\scriptscriptstyle_{|_P}}&=&0, \hbox{ R-P interface}.
\label{CR1}
\end{eqnarray}
Here ``H-P'' interface denotes the interface separating host from the plasmonic coating and ``R-P'' interface denotoes the interface separating the rod core material and the plasmonic coating and $n$ denotes
the normal vectors pointing from the core into the coating on the  ``R-P'' interface and the coating into the host on the ``H-P'' interface.
Expanding $\psi_1$ in terms of the complete set of orthonormal eigenfunctions $\{\psi_{\lambda_n}\}$ we obtain the  the representation 
\begin{eqnarray}
\psi_1=-\sum_{-1/2<\lambda_n<1/2}\left(\frac{(\alpha^1_{\lambda_n}+\epsilon_P^{-1}(\xi_0)\alpha^2_{\lambda_n})}{1+(\epsilon_P^{-1}(\xi_0)-1)(1-\lambda_n)}\right)\psi_{\lambda_n},\hbox{ in $Y\setminus R$}
\label{expansionpsi1}
\end{eqnarray}
with
\begin{eqnarray}
 \alpha^1_{\lambda_n}=\hat{\kappa}\cdot\int_H\nabla\psi_{\lambda_n}\,d\y, & \hbox{ and } \alpha^2_{\lambda_n}=\hat{\kappa}\cdot\int_P\nabla\psi_{\lambda_n}\,d\y.
\label{coefficients}
\end{eqnarray}
A straight forward calculation gives $\psi_0$ in $R$ in terms of the complete set of  Dirichlet eigenfunctions and eigenvalues $\{\mu_n\}$ and $\{\phi_n\}$:
\begin{eqnarray}
 &&\psi_0=\sum_{n=1}^{\infty}\frac{\mu_n<\phi_n>_R}{\mu_n-\xi_0}\phi_n , \hbox{ in $R$, with}
\label{explicitformofpsi0}
\end{eqnarray}
\begin{eqnarray}
 <\phi_n>_R=\int_R\phi_n\ d\y.
\label{coefficientsdir}
\end{eqnarray}

Setting  $v=1$ and $m=2$ in \eqref{summation2} we recover the solvability condition given by 
\begin{eqnarray}
&& \tau^2\int_{H\cup P}\big[ -\hat\kappa\cdot\xi_0\nabla \psi_1+\xi_0 \big]-\tau^2\epsilon_r\frac{\omega_p^2}{c^2}\int_H(-\hat\kappa\nabla \psi_1+1)\label{solvablety}
\\&&=\int_Y(\xi_0^2\psi_0-\epsilon_r\frac{\omega_p^2}{c^2}\xi_0 \psi_0)\nonumber
\end{eqnarray}
Substitution of the spectral representations for $\psi_1$ and $\psi_0$ given by \eqref{explicitformofpsi0} and \eqref{expansionpsi1} into \eqref{solvablety} delivers the homogenized dispersion relation
\begin{eqnarray}
 \xi_0=\tau^2 n_{eff}^{-2}(\xi_0),
\label{subwavelengthdispersion1nondim}
\end{eqnarray}
where the effective index of diffraction $n_{eff}^2$ depends upon the direction of propagation $\hat{\kappa}$ and is written
\begin{eqnarray}
n_{eff}^2(\xi_0)=\mu_{eff}(\xi_0)/\epsilon_{eff}^{-1}(\xi_0)\hat\kappa\cdot\hat\kappa.
\label{indexofrefract}
\end{eqnarray}
The frequency dependent effective magnetic permeability $\mu_{eff}$ and effective dielectric permittivity $\epsilon_{eff}$ are given by
\begin{eqnarray}
\mu_{eff}(\xi_0)=\int_Y\psi_0=\theta_H+\theta_P+\sum_{n=1}^\infty\frac{\mu_n<\phi_n>^2_R}{\mu_n-\xi_0}\label{effperm}
\end{eqnarray}
\par
and
\begin{eqnarray}
  \epsilon^{-1}_{eff}(\xi_0)\hat\kappa\cdot\hat\kappa&=& \int_{Y\setminus R}a_d^{-1}(y)(\nabla\psi_1+\hat{\kappa})\cdot\hat{\kappa}\,dy=\int_{Y\setminus R}a_d^{-1}(y)(\nabla\psi_1+\hat{\kappa})\cdot\overline{(\nabla\psi_1+\hat{\kappa})}\,dy\nonumber\\
 &=&\theta_H+\frac{\xi_0}{\xi_0-\frac{\epsilon_r\omega_p^2}{c^2}}\theta_P\label{effdielectricconst}\\
&&-\sum_{-1/2<\lambda_h <1/2}\left( \frac{\left(\xi_0-\frac{\epsilon_r\omega_p^2}{c^2} \right)|\alpha_{\lambda_h}^{(1)}|^2+2\frac{\epsilon_r\omega_p^2}{c^2}\alpha_{\lambda _h}^{(1)}\alpha_{\lambda_h}^{(2)}+\frac{\left(\frac{\epsilon_r\omega_p^2}{c^2}\right)^2}{\xi_0-\frac{\epsilon_r\omega_p^2}{c^2}}|\alpha_{\lambda_h}^{(2)}|^2}{\xi_0-(\lambda_h+\frac{1}{2})\frac{\epsilon_r\omega_p^2}{c^2}} \right),\nonumber
\end{eqnarray}
where $\theta_H$ and $\theta_P$ are the areas occupied by regions $H$ and $P$ respectively.

Writing out the dispersion relation \eqref{subwavelengthdispersion1nondim} explicitly in terms of $\mu_{eff}$ and $\epsilon_{eff}$ gives
\begin{eqnarray}
\mu_{eff}(\xi_0)\xi_0=\tau^2\epsilon_{eff}^{-1}(\xi_0)\hat{\kappa}\cdot\hat{\kappa}.
\label{longformnondindispersion}
\end{eqnarray}
From \eqref{longformnondindispersion} it is evident that there is a solution $\xi_0$ over intervals for which $\mu_{eff}(\xi_0)$ and $\epsilon_{eff}^{-1}(\xi_0)\hat{\kappa}\cdot\hat{\kappa}$ have the same sign.
It is also clear that there are an infinite number of intervals of the dispersion relation for which this is true and the branches of solutions to \eqref{longformnondindispersion} are  labeled $\{\xi_0^n\}_{n=1}^\infty$. We explicitly note the dependence of these branches on the wave number $k$ and propagation direction $\hat{\kappa}$ and write $\xi_0^n(k,\hat{\kappa})$. 
The power series for each branch of the dispersion relation and associated transverse magnetic Bloch wave solution of \eqref{Helmholtzr2} is given by
\begin{eqnarray}
 \xi^n=\xi_0^n(\tau,\hat{\kappa})+\sum_{l=1}^{\infty}(\tau\rho)^{l}\xi_l^n,
 \label{xipowerseriesnondimensional}
\end{eqnarray}
for
\begin{eqnarray}
\{-2\pi\leq \tau\rho\hat{\kappa}_1\leq 2\pi , -2\pi\leq \tau\rho\hat{\kappa}_2\leq 2\pi\}
\label{brullionfortaunondim}
\end{eqnarray}
and
\begin{eqnarray}
H_3^n=\underline u_0\left(\psi_0^n(\x/d)+\sum_{l=1}^{\infty}(\tau\rho)^li^l\psi_l^n(\x/d)\right)\exp\left\{i\left(k\hat\kappa\cdot \x-t\frac{\omega}{c}\right)\right\},
\label{summablenondim}
\end{eqnarray}
where
\begin{eqnarray}
\frac{\omega}{c}=\sqrt{\frac{\xi^n}{\epsilon_r}}.
\label{dispersionexplicitnondim}
\end{eqnarray}
For each branch of the dispersion relation,  the series converge for $\tau\rho$ sufficiently small, this follows from the theory developed in \cite{chenlipton}.

\section{Homogenization and energy flow for double negative effective properties}
The power series representation is used to show that branches of solutions \eqref{longformnondindispersion}  corresponding to $\mu_{eff}(\xi_0)<0$, $\epsilon_{eff}^{-1}(\xi_0)\hat{\kappa}\cdot\hat{\kappa}<0$ correspond to frequency intervals where the phase velocity in the effective medium is opposite to the direction of energy flow.
For H-polarized Bloch waves, the magnetic field $\textbf{H}(\x/d)=(0,0,H_3(\x /d))$ where $H_3(\x/d)$ is given by (\ref{summablenondim}) and the electric field $\textbf{E}(\x/d)=(E_1(\x/d), E_2(\x/d), 0)$. Both fields are related through (\ref{relationofefandmf}). Therefore 
\begin{eqnarray}
 \textbf{E}(\x/d)=\frac{ic}{\omega a_d}\partial_{x_2}H_3(\x/d){\bf e_1}-\frac{ic}{\omega a_d}\partial_{x_1}H_3(\x/d){\bf e_2},
\label{efield}
\end{eqnarray}
where $\bf {e_i}$ is the unit vector along the $x_i$ direction for $i=1,2,3$.
The time average of the Poynting vector is given by
\begin{eqnarray}
\textbf{P}^d&=&\frac{1}{2}Re[\textbf{E}(\x/d)\times \overline{\textbf{H}(\x/d)}]\nonumber\\
&=&\frac{1}{2}Re[E_2(\x/d)\overline{H_3(\x/d)}{\bf e_1}-E_1(\x/d)H_3(\x/d){\bf e_2}].
\label{poyntingvector}
\end{eqnarray}
Consider any fixed averaging domain $D$ transverse to the cylinders and the spatial average of the electromagnetic energy flow along the direction $\hat\kappa$  over this domain is written $\langle \textbf{P}\cdot\hat\kappa\rangle_D$.
Substituting (\ref{summablenondim}) and (\ref{efield}) into (\ref{poyntingvector}) and taking the limit of (\ref{poyntingvector}) as $d\rightarrow 0$ shows that the average electromagnetic energy flow along the direction $\hat\kappa$ is given by 
\begin{eqnarray}
\langle \textbf{P}\cdot\hat\kappa\rangle_D= \frac{1}{2}|\underline{u_0}|^2 n_{eff}\epsilon^{-1}_{eff}\hat\kappa\cdot\hat\kappa.
\label{poyntingalongk}
\end{eqnarray}
In the $d\rightarrow 0$ limit, the phase velocity is along the direction $\hat\kappa$ and determined by 
\begin{eqnarray}
\textbf v_p = \frac{c}{n_{eff}}\hat\kappa.
\label{phasevelocity}
\end{eqnarray}
Recall that we have pass bands over frequency intervals where $\epsilon^{-1}_{eff}\hat\kappa\cdot\hat\kappa$ and $\mu_{eff}$ are of the same sign. With this in mind equations (\ref{poyntingalongk}) and (\ref{phasevelocity}) show that in the homogenization limit the energy flow and phase velocity are in opposite directions over frequency intervals where the double negative property happens, i.e., $\epsilon^{-1}_{eff}\hat\kappa\cdot\hat\kappa<0$ and $\mu_{eff}<0$. These results are indicative of negative index behavior in the homogenization limit.

\section{Generalized electrostatic resonances for circular coated cylinders and the Rayleigh identity}
\label{Rayleigh's Identity}
In this section we develop a Rayleigh's identity for the eigenfunctions associated with generalized electrostatic resonances for coated circular cylinders. The Rayleigh method \cite{Rayleigh} has been generalized and applied to the analysis of wave propagation and the effective transport properties of composites and the approach taken here is motivated by the recent work \cite{McPhedranBook}, \cite{N} and \cite{Perrins}.
The electrostatic resonances $\lambda_h$ in (\ref{effdielectricconst}) are found by solving the following problem for the  potential $u$ inside a unit cell, i.e., $d=1$:
\begin{equation}
\begin{cases}
\Delta u=0 ~~~\text{ in  } H ,\\
\Delta u=0 ~~~\text{ in  } P ,\\
\end{cases}
\end{equation}
with the boundary conditions
\begin{eqnarray}
 \begin{cases}
  u|^{-}=u|^{+} ~~~\text{ on } \partial P,\\
\partial_{r}u|_{r=a}=0 ~~~\text{ on } \partial R,\\
\lambda[\partial_r u]^{-}_{+}=-\frac{1}{2}(\partial_ru^-+\partial_ru^+) \text{ on } \partial P ,\\
u\text{ is $Y$-periodic } . \\
 \end{cases}
\label{boundarycondition2}
\end{eqnarray}
Consequently, in polar coordinates $(r,\theta)$ , the expansions of the potential $u(r,\theta)$ are
\begin{eqnarray}
 u_p(r,\theta)=\sum_{l=1}^{\infty}(A_lr^l+B_lr^{-l})\cos l\theta  ~~~~\text{ in }  P ,
\label{solutioninP}
\end{eqnarray}
\begin{eqnarray}
 u_h(r,\theta)=\sum_{l=1}^{\infty}(C_lr^l+D_lr^{-l})\cos l\theta  ~~~~\text{ in }  H .
 \label{solutioninH}
\end{eqnarray}
\par
From the boundary conditions (\ref{boundarycondition2}), we can express $A_l, C_l $ and $D_l$ in terms of $B_l$.  Therefore we get
\begin{eqnarray}
 \begin{cases}
  A_l=a^{-2l}B_l ,\\
  C_l=\big(\frac{a^{2l}b^{-2l}-2\lambda}{1-2\lambda}\big) a^{-2l}B_l , \\
D_l=\big(\frac{a^{-2l}b^{2l}-2\lambda}{1-2\lambda}\big) B_l .\\
 \end{cases}
\label{intermsofB}
\end{eqnarray}
\par
The surface charge density $Q_s(\theta)$ is defined by
\begin{eqnarray}
 Q_s(\theta) &=& (\partial_ru_p-\partial_ru_h)|_{r=b} \nonumber 
\\&=& 2\sum_{l=1}^{\infty}\big(\frac{b^{l-1}a^{-2l}-b^{-(l+1)}}{1-2\lambda}\big)lB_l\cos l\theta .
\label{density}
\end{eqnarray}
The potential at an arbitary point $G(r,\theta) $ is given by
\begin{eqnarray}
 u(r,\theta)=-\frac{1}{2\pi}\sum_{j}\int_{\partial s_j} Q_{s_j}^j(\bf{t_{(s_j)}})\ln (|\bf{r}-\bf{t_{(s_j)}}|)d s_j .
\label{potential}
\end{eqnarray}

In the summation $j$ refers to the $j$th cylinder and the vector $\bf{t_{(s_j)}}$ extended from the origin to the area element $ds_j$ on its shell. $Q_{s_j}^j(\bf{t_{(s_j)}})$ is the surface charge density at this area element. We sum over all cylinders in the lattice and integrate over the entire surface of each cylinder. We introduce the vectors $\boldsymbol{\rho_j}$ pointing from the center of the $j$th cylinder to the point $G$ and $\bf{s}$ pointing from the center of the $j$th cylinder to $\partial{s_j}$ . 
\begin{figure}[h!]
\begin{center}
\begin{psfrags} 
\psfrag{x}{$x$} 
\psfrag{y}{$y$}
\psfrag{Rj}{$\bf{R_j}$}
\psfrag{o}{$O$}
\psfrag{G}{$G$}
\psfrag{S}{$S$}
\psfrag{phi}{$\phi_j$}
\psfrag{thetaj}{$\theta_j$}
\psfrag{rho}{$\boldsymbol{\rho_j}$}
\psfrag{s}{$\bf{s}$}
\psfrag{t}{$\bf{t_{(s_j)}}$}
\psfrag{r}{$\bf{r}$}
\includegraphics[width=0.6\textwidth]{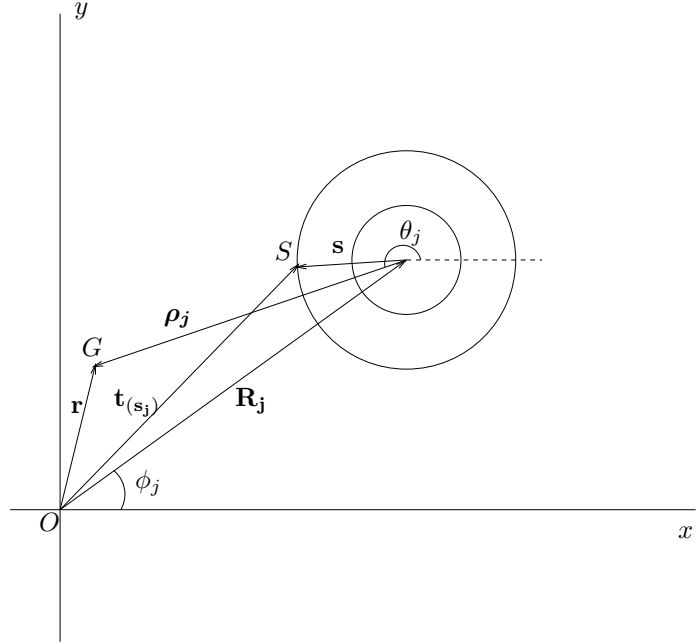}
\end{psfrags} 
\caption{$G$ is a typical field point while $S$ present a point on the shell boundary of the $j$th cylinder.}
\label{geometry}
\end{center}
\end{figure}

From Figure \ref{geometry} , it is easy to see that
\begin{eqnarray}
 \bf{r}-\bf{t_{(s_j)}}=\boldsymbol{\rho_j}-\bf{s} . 
\end{eqnarray}
\par
First suppose that the field point $G$ is in the shell of the central cylinder ,i.e., $a<r<b$.  Since $|\bf{s}|=b , \bf{\rho_0}=\bf{r}$  and  $|\boldsymbol{\rho_j}|>b~~\forall j\neq 0 , ~a<|\bf{r}|<b~~\text{for } j=0$.  Then we can expand the logarithms in (\ref{potential}) 
\begin{eqnarray}
 \ln |\bf{r}-\bf{t_{(s_j)}}|&=&\ln |\boldsymbol{\rho_j}-\bf{s}|\nonumber
\\& =&
\begin{cases}
 \ln \rho_j-\sum_{n=1}^{\infty}\frac{1}{n}(\frac{b}{\rho_i})^n\cos n(\theta'-\theta_j)~~~ j\neq 0,\\
 \ln b-\sum_{n=1}^{\infty}\frac{1}{n}(\frac{r}{b})^n\cos n(\theta'-\theta)~~~~~~ j=0 ,\\
\end{cases}
\label{expansionoflog}
\end{eqnarray}
where $\theta'$ is the polar angle defining the orientation of $\bf{s}$ and $\theta_j$ specifies the direction of $\boldsymbol{\rho_j}$ . For the $j$th cylinder the surface charge density takes the form (\ref{density}) depending on $\theta'$. From (\ref{potential}), (\ref{expansionoflog}) and the orthonormality properties of Sines and Cosines we find ($ds_j=bd\theta' $)
\begin{eqnarray}
 u(r,\theta)=\sum_{l=1}^{\infty}\frac{1}{1-2\lambda}(a^{-2l}-b^{-2l})B_lr^l\cos l\theta 
+\sum_{j\neq 0}\sum_{l=1}^{\infty}\frac{1}{1-2\lambda}\big(\frac{a^{-2l}b^{2l}-1}{\rho_j^l}\big)B_l \cos l\theta_j 
\label{summationofpotential}
\end{eqnarray}
for $a<r<b$.
\par
In (\ref{summationofpotential}) the first sum is from the central cylinder and the second sum is over all the other cylinders.  In the shell , we must have
\begin{eqnarray}
 u(r,\theta)=u_p(r,\theta)=\sum_{l=1}^{\infty}(A_lr^l+B_lr^{-l})\cos l\theta .
\label{solutioninshell}
\end{eqnarray}
Plugging (\ref{intermsofB}) and (\ref{summationofpotential}) into (\ref{solutioninshell}) , we obtain the Rayleigh's identity for a square array of coated cylinders:
\begin{eqnarray}
 \sum_{l=1}^{\infty}\left [\left(\frac{b^{-2l}-2\lambda a^{-2l}}{1-2\lambda}\right)r^l+r^{-l}\right]B_l \cos l\theta = \sum_{j\neq 0}\sum_{l=1}^{\infty}\left(\frac{a^{-2l}b^{2l}-1}{1-2\lambda}\right)\frac{B_l}{\rho_j^l}\cos l\theta_j
\label{rayleighidentity}
\end{eqnarray}
for $a<r<b$.
\par
Next we suppose that the field point $G$ is outside the shell of the the central cylinder, i.e., $r>b$. The proof is similar to the case in the shell $P$. Since $r>b$ ,  (\ref{expansionoflog})  changes to be
\begin{eqnarray}
 \ln |\bf{r}-\bf{t_{(s_j)}}|&=&\ln |\boldsymbol{\rho_j}-\bf{s}|\nonumber
\\& =&
\begin{cases}
 \ln \rho_j-\sum_{n=1}^{\infty}\frac{1}{n}(\frac{b}{\rho_i})^n\cos n(\theta'-\theta_j)~~~ j\neq 0,\\
 \ln r-\sum_{n=1}^{\infty}\frac{1}{n}(\frac{b}{r})^n\cos n(\theta'-\theta)~~~~~~ j=0 .\\
\end{cases}
\label{expansionoflog2}
\end{eqnarray}
Then like (\ref{summationofpotential})  , we have
\begin{eqnarray}
 u(r,\theta)=\sum_{l=1}^{\infty}\frac{1}{1-2\lambda}(a^{-2l}b^{2l}-1)B_lr^{-l}\cos l\theta 
+\sum_{j\neq 0}\sum_{l=1}^{\infty}\frac{1}{1-2\lambda}\big(\frac{a^{-2l}b^{2l}-1}{\rho_j^l}\big)B_l \cos l\theta_j 
\label{summationofpotential2}
\end{eqnarray}
for $r>b$. In the host $H$, we must have
\begin{eqnarray}
 u(r,\theta)=u_h(r,\theta)=\sum_{l=1}^{\infty}(C_lr^l+D_lr^{-l})\cos l\theta .
\label{solutioninhost}
\end{eqnarray}
A  calculation similar to the previous in the shell shows that (\ref{rayleighidentity}) holds for $r>b$. Therefore (\ref{rayleighidentity}) is true for all $r>a$.

\section{Numerical calculation of generalized electrostatic spectra for coated cylinders}
\label{numericalresult}
Now we set $r=b$ in (\ref{rayleighidentity}) and apply Rayleigh's method \cite{Rayleigh} to recover a linear system of equations to determine $\lambda$ and $\{B_l\}_{l=1}^{\infty}$. Let $R_j$ be the distance from the origin to the center of the $j$th cylinder and the $\phi_j$ be the polar angle defining the orientation of $\bf{R_j}$ (See Figure \ref{geometry} ). Then for $r=b$ we have 
\begin{eqnarray}
 \rho_j\cos \theta_j+i\rho_j\sin \theta_j= (b\cos\theta+ib\sin\theta) -(R_j\cos\phi_j+iR_j\sin\phi_j)
\end{eqnarray}
The Rayleigh's identity (\ref{rayleighidentity}) with $r=b$ becomes
\begin{small}\begin{eqnarray}
 &&Re \sum_{l=1}^{\infty}\left [\left(\frac{b^{-2l}-2\lambda a^{-2l}}{1-2\lambda}\right)B_l(b\cos\theta+ib\sin\theta)^l+B_l(b\cos\theta+ib\sin\theta)^{-l}\right]\nonumber
\\&&=Re \sum_{j\neq 0}[\left(\frac{ a^{-2}b^2-1}{1-2\lambda}\right)B_1\left((b\cos\theta+ib\sin\theta) -(R_j\cos\phi_j+iR_j\sin\phi_j)\right)^{-1}\nonumber
\\&&~+\left(\frac{ a^{-4}b^4-1}{1-2\lambda}\right)B_2\left((b\cos\theta+ib\sin\theta) -(R_j\cos\phi_j+iR_j\sin\phi_j)\right)^{-2}\nonumber
\\&&~+\left(\frac{ a^{-8}b^8-1}{1-2\lambda}\right)B_3\left((b\cos\theta+ib\sin\theta) -(R_j\cos\phi_j+iR_j\sin\phi_j)\right)^{-3}\nonumber
\\&&~+\cdots]
\label{expandrayleigh}
\end{eqnarray}              \end{small}
Using the generalized binomial theorem on the right-hand side of (\ref{expandrayleigh}), we have
\begin{small}\begin{eqnarray}
 &&Re \sum_{l=1}^{\infty}\left [\left(\frac{b^{-2l}-2\lambda a^{-2l}}{1-2\lambda}\right)B_l(b\cos\theta+ib\sin\theta)^l+B_l(b\cos\theta+ib\sin\theta)^{-l}\right]\nonumber
\\&&=Re \sum_{j\neq 0}[\left(\frac{ a^{-2}b^2-1}{1-2\lambda}\right)B_1(-1)^{1}(R_j\cos\phi_j+iR_j\sin\phi_j)^{-1}\sum_{k=0}^{\infty}\left(\frac{b\cos\theta+ib\sin\theta}{R_j\cos\phi_j+iR_j\sin\phi_j}\right)^{k}\nonumber
\\&&~+\left(\frac{ a^{-4}b^4-1}{1-2\lambda}\right)B_2(-1)^{2}(R_j\cos\phi_j+iR_j\sin\phi_j)^{-2}\sum_{k=0}^{\infty}\binom {2+k-1} {k}\left(\frac{b\cos\theta+ib\sin\theta}{R_j\cos\phi_j+iR_j\sin\phi_j}\right)^{k}\nonumber
\\&&~+\left(\frac{ a^{-8}b^8-1}{1-2\lambda}\right)B_3(-1)^{3}(R_j\cos\phi_j+iR_j\sin\phi_j)^{-3}\sum_{k=0}^{\infty}\binom {3+k-1} {k}\left(\frac{b\cos\theta+ib\sin\theta}{R_j\cos\phi_j+iR_j\sin\phi_j}\right)^{k}\nonumber
\\&&~+\cdots].
\label{expandrayleigh2}
\end{eqnarray}              \end{small}
Equating the coefficients of $\cos l\theta$ in (\ref{expandrayleigh2}) between left- and right-hand sides, we obtain
\begin{eqnarray}
 &&\left(\frac{b^{-2l}-2\lambda a^{-2l}}{1-2\lambda}\right)B_lb^l+B_lb^{-l}\nonumber
\\&&=Re\{ \sum_{j\neq 0}[\left(\frac{ a^{-2}b^2-1}{1-2\lambda}\right)B_1(-1)^{1}(R_j\cos\phi_j+iR_j\sin\phi_j)^{-l-1}b^l\nonumber
\\&&~+\left(\frac{ a^{-4}b^4-1}{1-2\lambda}\right)B_2(-1)^{2}(R_j\cos\phi_j+iR_j\sin\phi_j)^{-l-2}\binom {2+l-1} {l}b^l\nonumber
\\&&~+\left(\frac{ a^{-8}b^8-1}{1-2\lambda}\right)B_3(-1)^{3}(R_j\cos\phi_j+iR_j\sin\phi_j)^{-l-3}\binom {3+l-1} {l}b^l\nonumber
\\&&~+\cdots]\}\nonumber
\\&&=\sum_{j\neq 0}\sum_{m=1}^{\infty}\left[\left(\frac{ a^{-2m}b^{2m}-1}{1-2\lambda}\right)\binom {m+l-1} {l}B_m(-1)^{m}\frac{\cos (l+m)\phi_j}{R_j^{l+m}}b^l\right]\nonumber
\\&&=\sum_{m=1}^{\infty}\left[\left(\frac{ a^{-2m}b^{2m}-1}{1-2\lambda}\right)\binom {m+l-1} {l}B_m(-1)^{m}S_{l+m}b^l\right],
\label{equatingcos}
\end{eqnarray}
where the quantities $S_n$ are the lattice sums
\begin{eqnarray}
 S_n=\sum_{j\neq 0}\frac{\cos n\phi_j}{R_j^n}.
\end{eqnarray}
A list of numerical values for  $S_n$ is tabulated in the paper of Perrins Perrins et. al. \cite{Perrins}. Rewriting (\ref{equatingcos}) as
\begin{eqnarray}
 \lambda B_l=\left(\frac{b^{-2l}}{a^{-2l}+b^{-2l}}\right)B_l-\frac{1}{2(a^{-2l}+b^{-2l})}\sum_{m=1}^{\infty}( a^{-2m}b^{2m}-1)\binom {m+l-1} {l}(-1)^{m}S_{l+m}B_m.
\label{systemequation}
\end{eqnarray}
The system (\ref{systemequation}) may be written  in the matrix form
\begin{eqnarray}
 AB=\lambda B,
\label{matrixform}
\end{eqnarray}
where $B=(B_1,B_2,B_3, \cdots)^{T}$ and the infinite dimentional matrix $A$ has the elements
\begin{equation}
 A_{lm}=
\begin{cases}
 \frac{b^{-2l}}{a^{-2l}+b^{-2l}}-\frac{1}{2(a^{-2l}+b^{-2l})}( a^{-2m}b^{2m}-1)\binom {m+l-1} {l}(-1)^{m}S_{l+m}~~~~~~  l=m,\\
-\frac{1}{2(a^{-2l}+b^{-2l})}( a^{-2m}b^{2m}-1)\binom {m+l-1} {l}(-1)^{m}S_{l+m} ~~~~~~~~~~~~~~~~~~~l\neq m.\\
\end{cases}
\end{equation}
\par
We solve (\ref{matrixform}) numerically,  truncating the sum after $N$ terms to  find the approximation of the potential $u$ and the generalized electrostatic resonances.  Table \ref{Table1} gives the eigenvalues corresponding to different $N$ for $a=0.2$ and $b=0.4$.
\begin{table}[h!]
\begin{center}
\begin{tabular}{|c|c|c|c|}
  \hline
   & $N=10$ & $N=15$ & $N=20$\\
  \hline
  &$3.5080\times 10^{-1}$ &$3.5080\times 10^{-1}$ &$3.5080\times 10^{-1}$ \\
  &$1.5379\times 10^{-2}$ &$1.5379\times 10^{-2}$ &$1.5379\times 10^{-2}$   \\
  &$9.7557\times 10^{-4}$ &$9.7557\times 10^{-4}$ &$9.7557\times 10^{-4}$  \\
  &$6.1031\times 10^{-5}$ &$6.1031\times 10^{-5}$  &$6.1031\times 10^{-5}$    \\
  &$3.8147\times 10^{-6}$  &$3.8147\times 10^{-6}$ &$3.8147\times 10^{-6}$  \\
  &$-2.0285\times 10^{-3}$ &$2.3842\times 10^{-7}$ &$2.3842\times 10^{-7}$  \\
  &$-5.5339\times 10^{-3}$ &$1.4901\times 10^{-8}$ &$1.4901\times 10^{-8}$  \\
  &$-1.5014\times 10^{-2}$ &$9.3132\times 10^{-10}$ &$9.3132\times 10^{-10}$ \\
  &$-4.4538\times 10^{-2}$ &$-2.8905\times 10^{-10}$ & $5.8208\times 10^{-10}$   \\
  &$-4.7947\times 10^{-2}$ &$-7.6128\times 10^{-10}$ &$3.6380\times 10^{-10} $  \\
 $\lambda$ &                     &$-2.0285\times 10^{-3}$ & $-1.6665\times 10^{-5}$   \\
  &                     &$-5.5339\times 10^{-3}$ &$-4.2856\times 10^{-5}$   \\
  &                     &$-1.5014\times 10^{-2}$ &$-1.1088\times 10^{-4}$   \\
  &                     &$-4.4538\times 10^{-2}$ &$-2.8905\times 10^{-4}$   \\
  &                     &$-4.7947\times 10^{-2}$ &$-7.6128\times 10^{-4}$  \\
  &                     &                      &$-2.0285\times 10^{-3}$  \\
  &                     &                      &$-5.5339\times 10^{-3}$  \\
  &                     &                      &$-1.5014\times 10^{-2}$  \\
  &                     &                      &$-4.4538\times 10^{-2}$  \\
  &                     &                      &$-4.7947\times 10^{-2}$ \\
   \hline
\end{tabular}
\end{center}
\caption{The eigenvalues corresponding to $N=10, 15,20$ with $a=0.2$ , $b=0.4$.}
\label{Table1}
\end{table}
These numerical results confirm that the eigenvalues have an accumulation point at $0$. For illustration, if $\lambda$ is the first eigenvalue $3.5080\times 10^{-1}$ with $a=0.2$ and $b=0.4$ , then we notice that $B_1\approx 1$ and the remaining coefficients $B_l$'s are close to $0$. Hence (\ref{summationofpotential}) and (\ref{summationofpotential2}) show that the potential $u(r,\theta)$ is well approximated by
\begin{eqnarray}
u(r,\theta)\approx
\begin{cases}
 \frac{1}{1-2\lambda}(a^{-2}-b^{-2})r\cos \theta 
+\sum_{j\neq 0}\frac{1}{1-2\lambda}\big(\frac{a^{-2}b^{2}-1}{\rho_j}\big) \cos \theta_j ~~~~~~\text{in } P,\\
\frac{1}{1-2\lambda}(a^{-2}b^{2}-1)r^{-1}\cos \theta 
+\sum_{j\neq 0}\frac{1}{1-2\lambda}\big(\frac{a^{-2}b^{2}-1}{\rho_j}\big)\cos \theta_j ~~~\text{in } H.\\
\end{cases}
\end{eqnarray}
Equivalently (\ref{solutioninP}) and (\ref{solutioninH}) show
\begin{eqnarray}
u(r,\theta)\approx
\begin{cases}
 a^{-2}r\cos \theta 
+r^{-1}\cos\theta ~~~~~~~~~~~~~~~~~~~~~~~~~~~~~~~~\text{in } P,\\
\big(\frac{a^{2}b^{-2}-2\lambda}{1-2\lambda}\big) a^{-2} r\cos \theta 
+\big(\frac{a^{-2}b^{2}-2\lambda}{1-2\lambda}\big) r^{-1}\cos\theta ~~~\text{in } H.\\
\end{cases}
\end{eqnarray}
\par
\begin{figure}[h!]
\centering
\mbox{\subfigure[]{\epsfig{figure=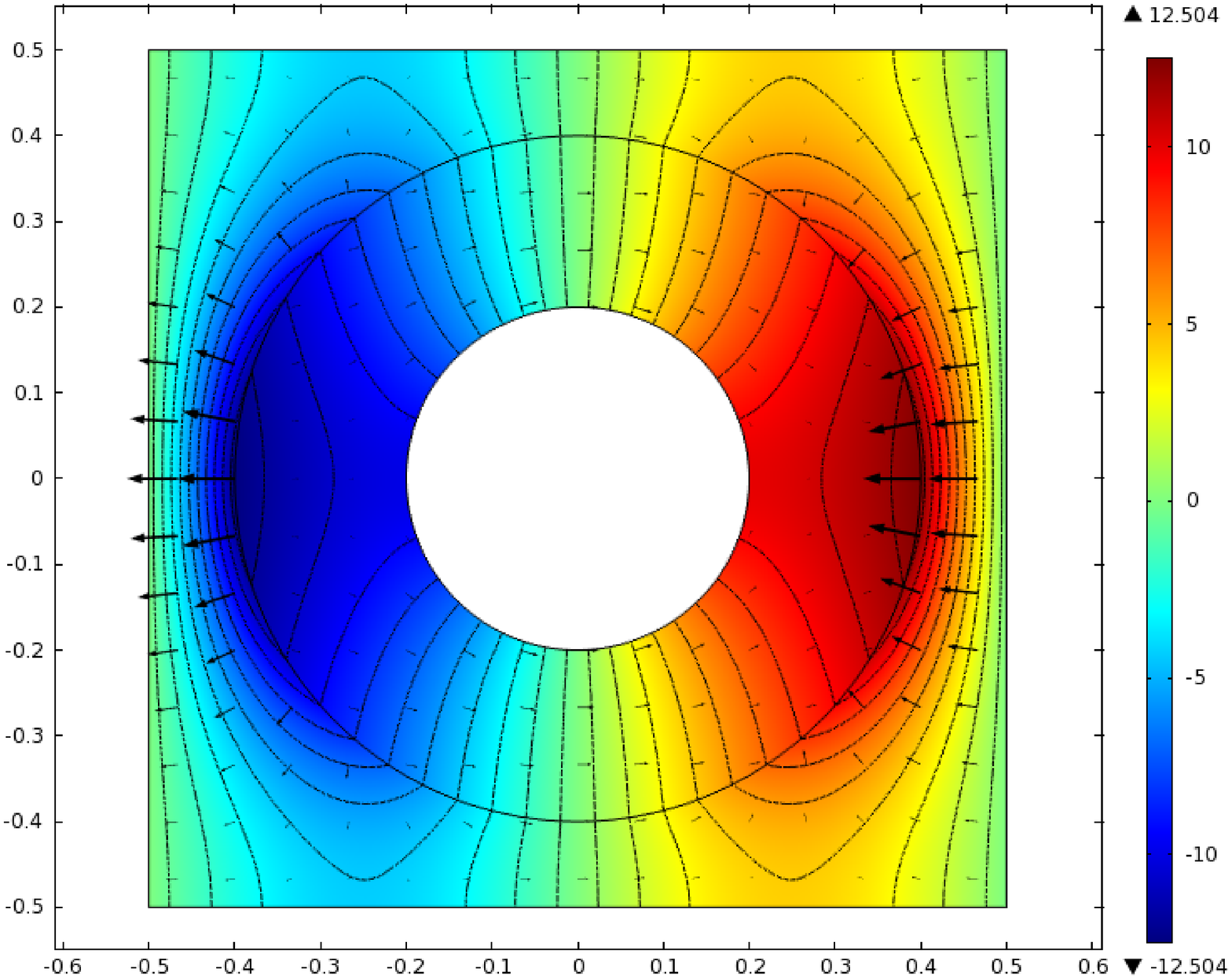 , width=2.5in}}\quad\quad\quad\quad\quad
\subfigure[]{\epsfig{figure=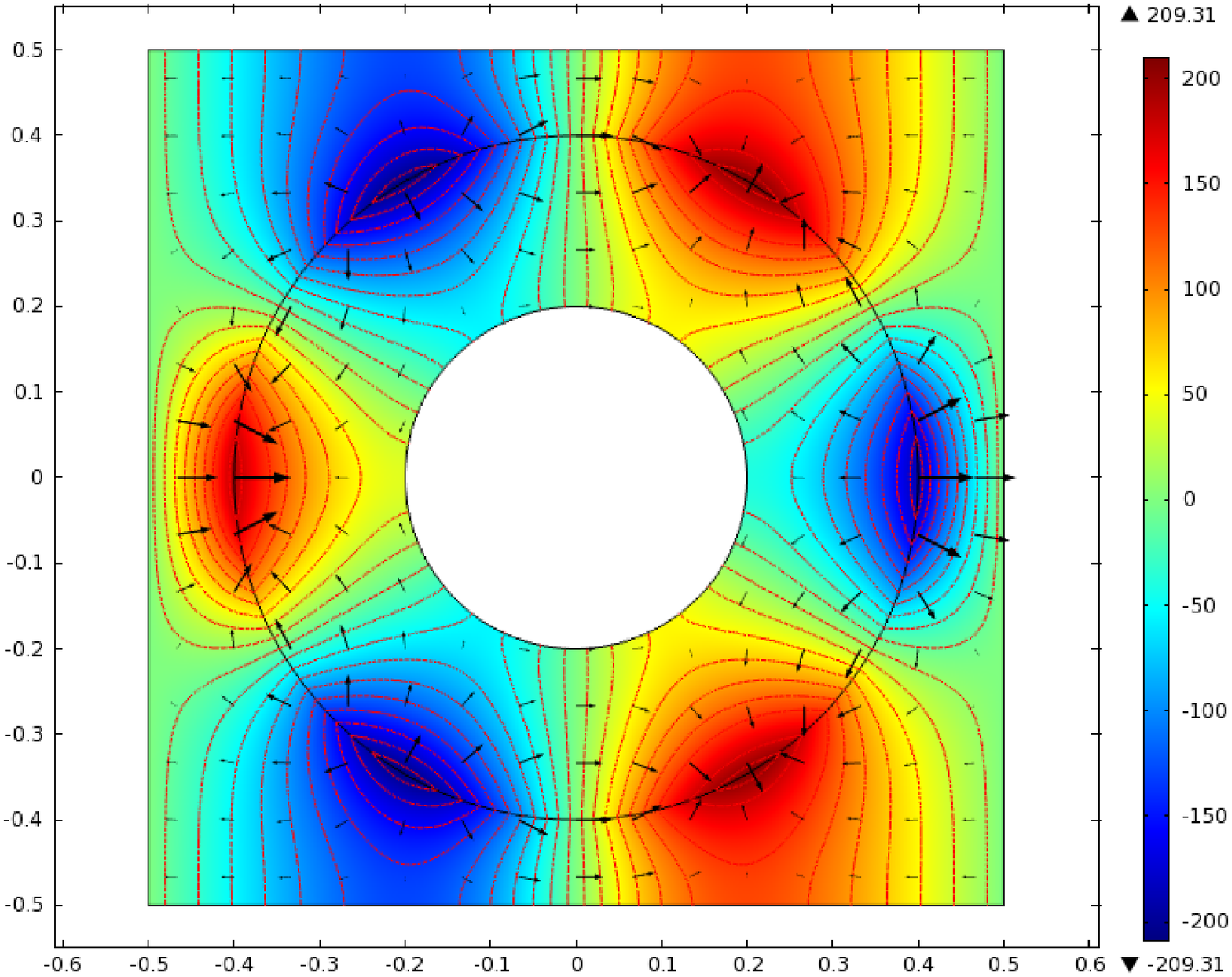 , width=2.5in}}}
\caption{(a) the solution corresponding to the eigenvalue $\lambda=3.5080\times 10^{-1}$; (b) the solution corresponding to the eigenvalue $\lambda=1.5379\times 10^{-2}$.}
\label{eigenvaluesolution}
\end{figure}

The the solution $u$ corresponding to the first two eigenvalues with $a=0.2$ and $b=0.4$ are illustrated in Figure \ref{eigenvaluesolution}.

\par
\section{Numerical calculation of the dispersion relation and comparison with power series}
\label{simulation}
In this section we verify that the leading order dispersion relation expressed in terms of effective properties is a good predictor of the dispersive behavior of the metamaterial for periods with finite size $d>0$. The usefulness of the effective properties for predicting metamaterial behavior away from the homogenization limit can be explicitly seen from power series formula for the dispersion relation.
To proceed we fix  $d=c/\omega_p$ and the dimensionless ratio $\rho=d/\sqrt{\epsilon_r}$. With this choice of variables the frequency dependent effective magnetic permeability $\mu_{eff}$ and effective dielectric permittivity $\epsilon_{eff}$ are written as
\begin{eqnarray}
\mu_{eff}(\omega_0/\omega_p)=\int_Y\psi_0=\theta_H+\theta_P+\sum_{n=1}^\infty\frac{\mu_n<\phi_n>^2_R}{\rho^{-2}\left(\mu_n\rho^2-(\frac{\omega_0}{\omega_p})^2\right)}\label{effperm2}
\end{eqnarray}
and
\begin{eqnarray}
  \epsilon^{-1}_{eff}(\omega_0/\omega_p)\hat\kappa\cdot\hat\kappa &=&\theta_H+\frac{(\frac{\omega_0}{\omega_p})^2}{(\frac{\omega_0}{\omega_p})^2-1}\theta_P\label{effdielectricconst2}\\
&&-\sum_{-\frac{1}{2}<\lambda_h <\frac{1}{2}}\left( \frac{\left((\frac{\omega_0}{\omega_p})^2-1 \right)^2|\alpha_{\lambda_h}^{(1)}|^2+2\left((\frac{\omega_0}{\omega_p})^2-1\right)\alpha_{\lambda _h}^{(1)}\alpha_{\lambda_h}^{(2)}+|\alpha_{\lambda_h}^{(2)}|^2}{\left((\frac{\omega_0}{\omega_p})^2-(\lambda_h+\frac{1}{2})\right)\left((\frac{\omega_0}{\omega_p})^2-1\right)} \right).\nonumber
\end{eqnarray}
In these variables the leading order dispersion relation is given by 
\begin{eqnarray}
 (dk)^2=(\frac{\omega_0}{\omega_p})^2n_{eff}^{2},
\label{subwavelengthdispersion1}
\end{eqnarray}
where the effective index of diffraction $n_{eff}^2$ depends upon the direction of propagation $\hat{\kappa}$ and normalized frequency $\frac{\omega_0}{\omega_p}$ and is written
\begin{eqnarray}
n_{eff}^2=\mu_{eff}(\frac{\omega_0}{\omega_p})/(\epsilon_{eff}^{-1}(\frac{\omega_0}{\omega_p})\hat\kappa\cdot\hat\kappa).
\label{indexofrefract2}
\end{eqnarray}
The dispersion relation for the metamaterial crystal in the new variables is given by 
\begin{eqnarray}
 \left(\frac{\omega}{\omega_p}\right)^2=\left(\frac{\omega_0}{\omega_p}\right)^2+\sum_{l=1}^{\infty}(dk)^{l}\left(\frac{\omega_l}{\omega_p}\right)^2,
 \label{xipowerseriesnew}
\end{eqnarray}
for
\begin{eqnarray}
\{-2\pi\leq dk\hat{\kappa}_1\leq 2\pi , -2\pi\leq dk\hat{\kappa}_2\leq 2\pi\}.
\label{brullionfortau}
\end{eqnarray}
It is clear from \eqref{xipowerseriesnew} that the roots $\frac{\omega_0}{\omega_p}$ of the effective dispersion relation \eqref{subwavelengthdispersion1} determines the leading order dispersive behavior for periods of finite size. We point out that the notion of effective properties for metamaterials has proved to be an elegant concept for explaining experimental results. Here it is seen that the power series 
\eqref{xipowerseriesnew} exhibits the precise way in which effective properties influence leading order behavior for period cells of size $d>0$.

\begin{figure}[!htb]
\centering
\subfigure[]{\epsfig{figure=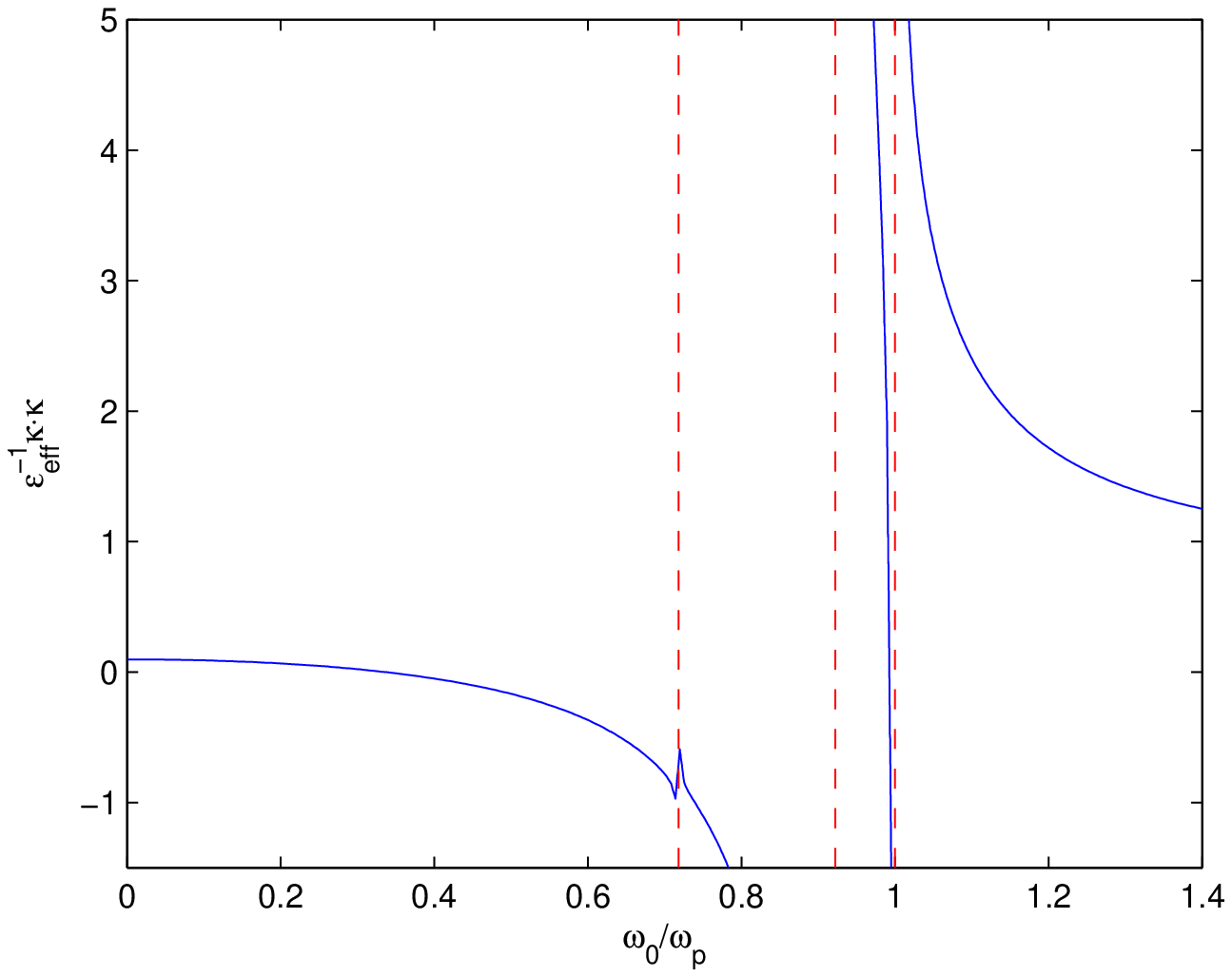 , width=1.9in}}~
\subfigure[]{\epsfig{figure=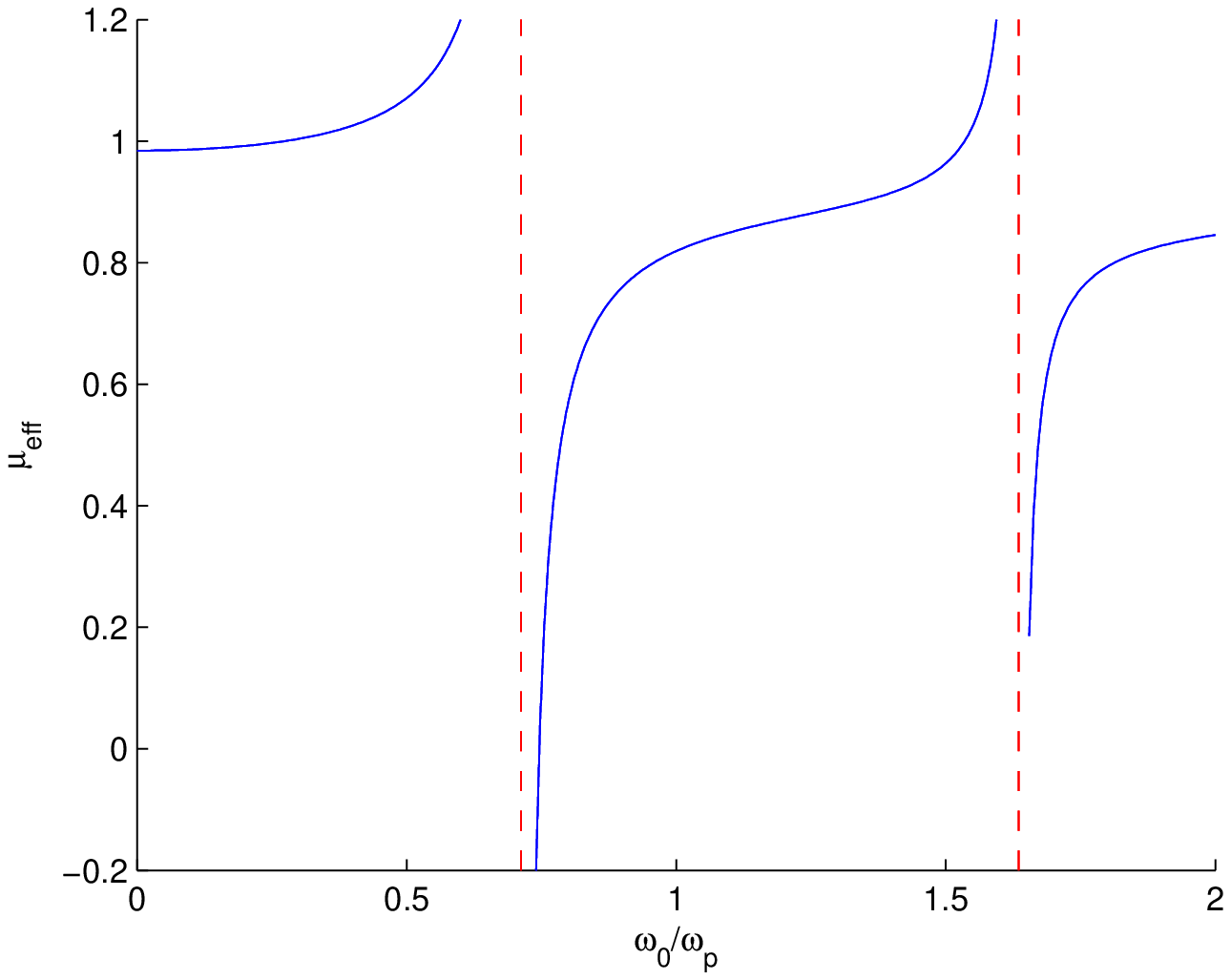 , width=1.9in}}~
\subfigure[]{\epsfig{figure=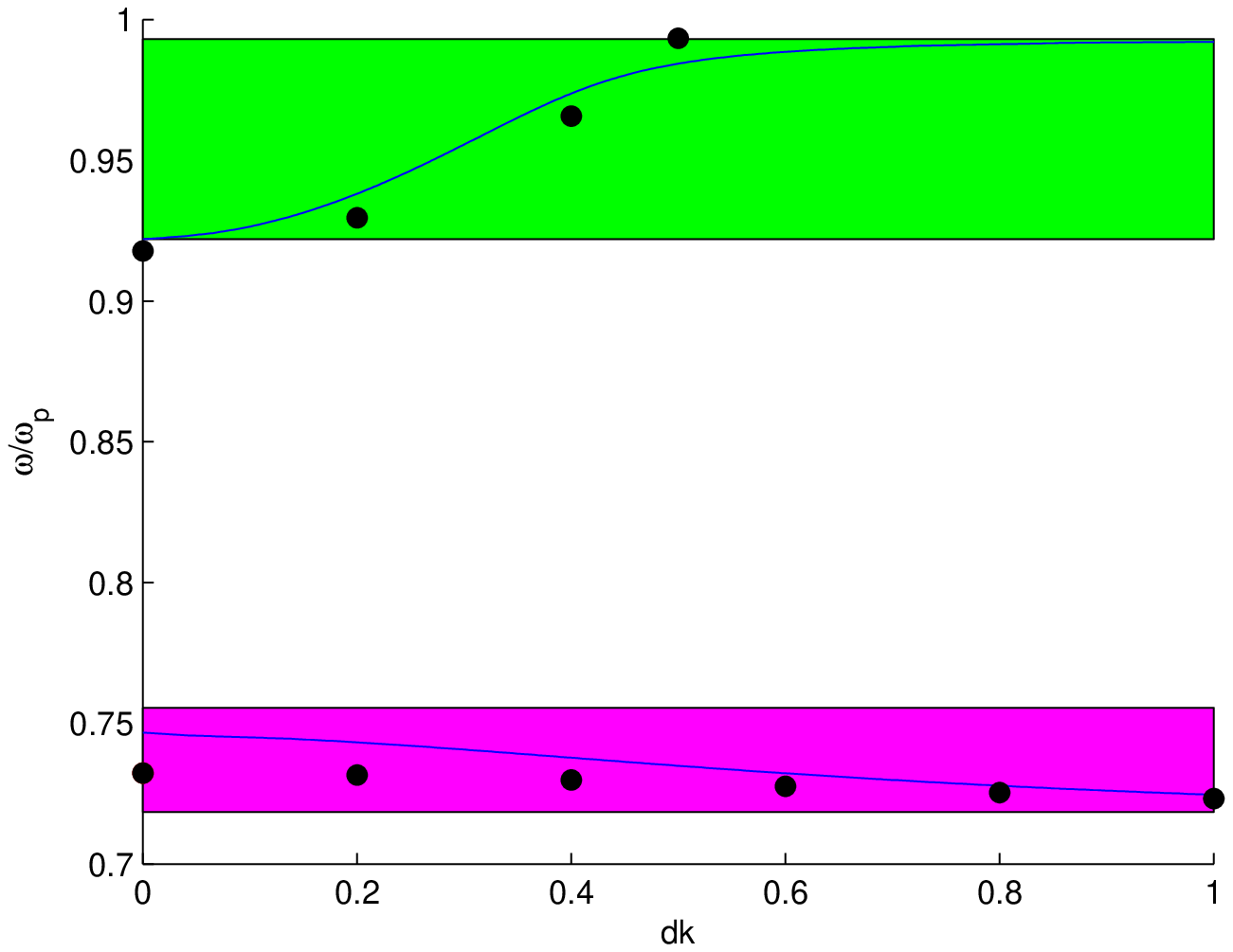 , width=1.9in}}
\caption{the case of $a=0.2d, b=0.4d$ and $\epsilon_R=285$.  
Notice that the vertical dash lines are the asymptotes.
}
\label{example1}
\end{figure}

\begin{figure}[!htb]
\centering
\subfigure[]{\epsfig{figure=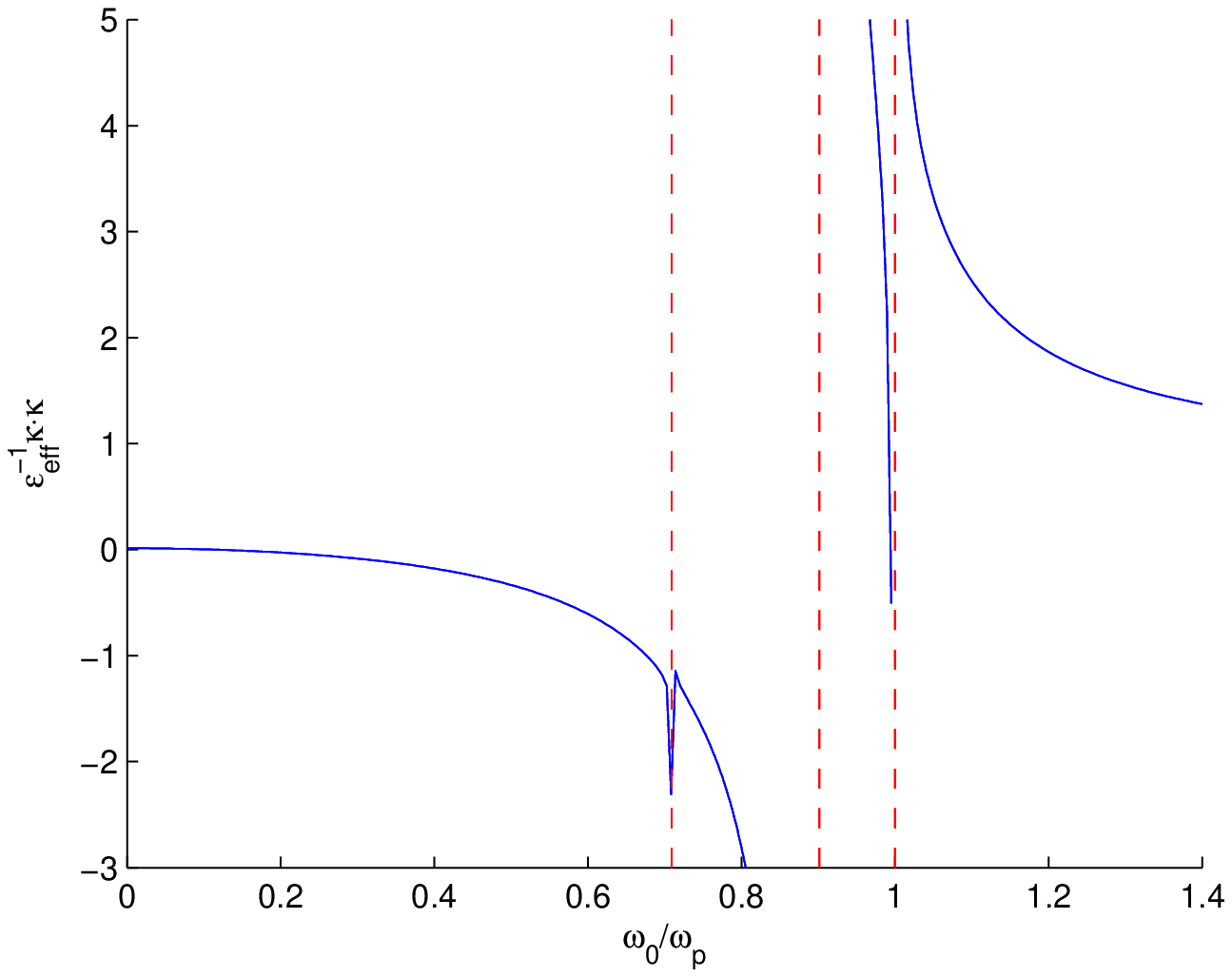 , width=1.9in}}~
\subfigure[]{\epsfig{figure=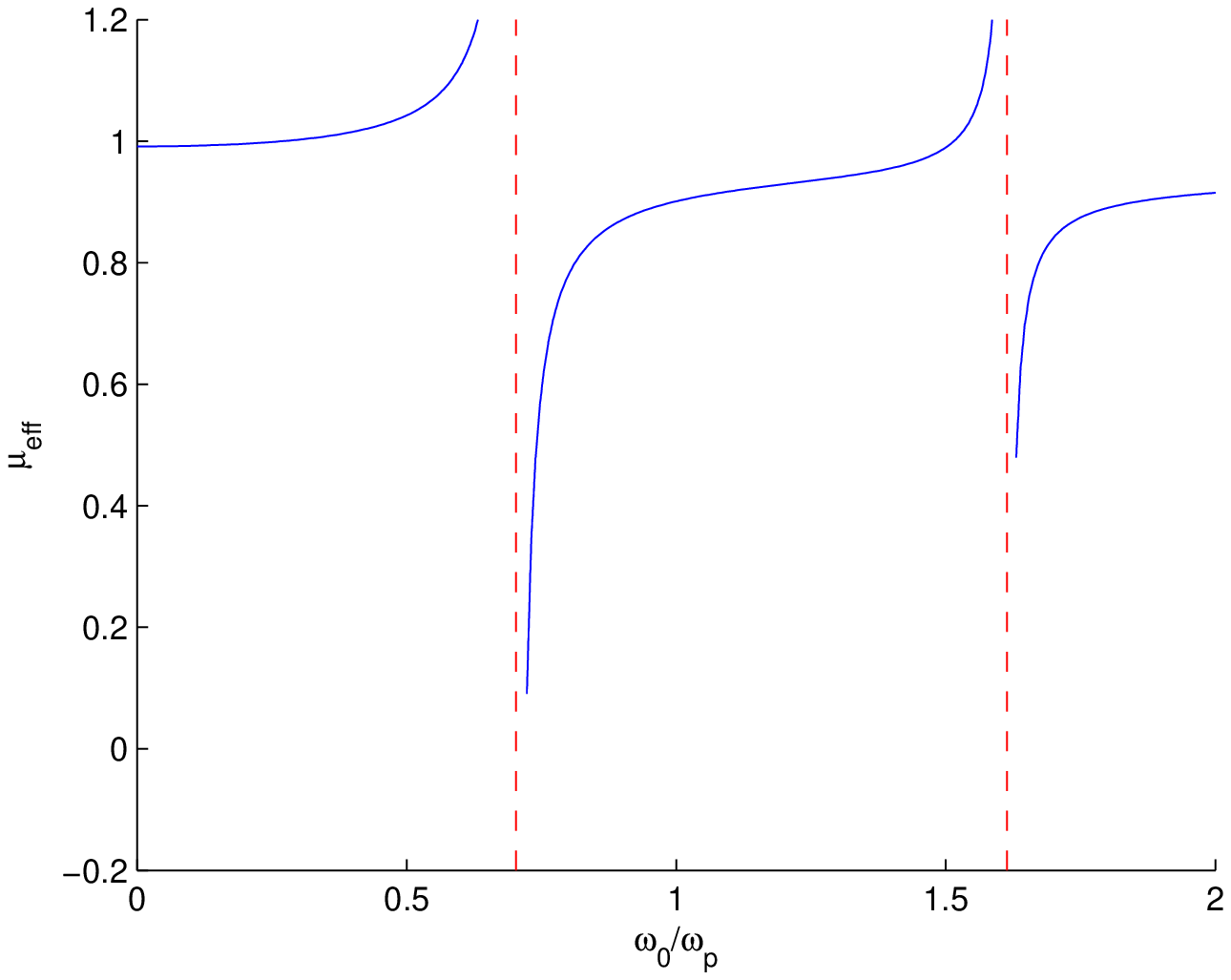 , width=1.9in}}~
\subfigure[]{\epsfig{figure=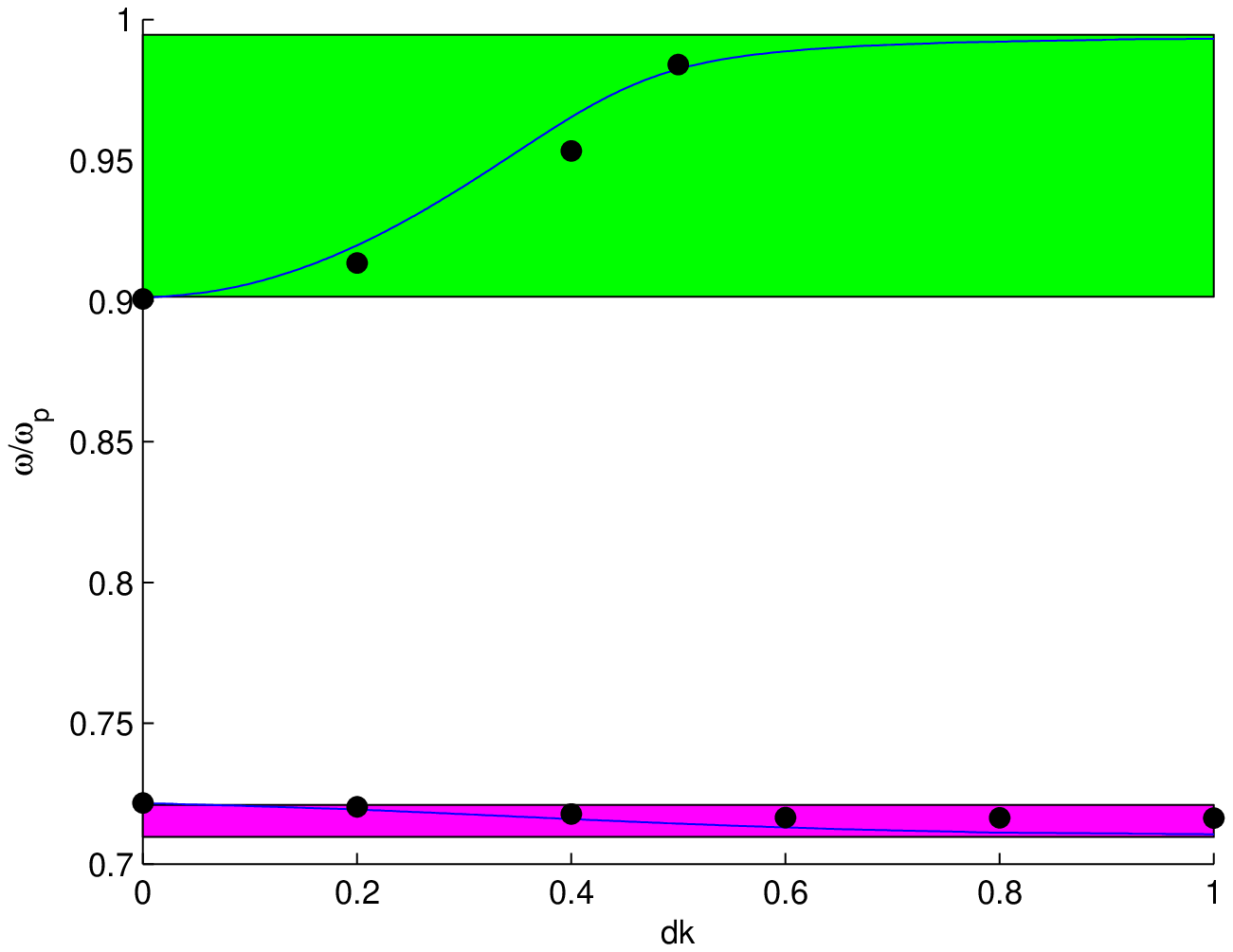 , width=1.9in}}
\caption{the case of $a=0.15d, b=0.4d$ and $\epsilon_R=285$. 
Notice that the vertical dash lines are the asymptotes.
}
\label{example2}
\end{figure}

To verify that the leading order dispersion relation is a good predictor of the dispersive behavior of the metamaterial,
we numerically compute the solutions $(\omega/\omega_p)^2$ and $u^d$  of the nonlinear eigenvalue problem given by equation (\ref{Helmholtzd}). These computations are carried out for different wave numbers $k$ with propagation along the direction $(1,0)$. The simulations are carried out using COMSOL software. Here we are interested in the range  $\omega/\omega_p<1$ for which the plasmonic material has a negative permittivity, $\epsilon_P<0$. 
Two examples are considered:
the first example is the case of  $a=0.2d, b=0.4d$ and $\epsilon_R=\epsilon_r/d^2=285$ and the second example is for $a=0.15d, b=0.4d$ and $\epsilon_R=285$. Figure \ref{example1} (a)(b) (Figure\ref{example2} (a)(b)) are the graphs of the effective properties $\epsilon^{-1}_{eff}(\frac{\omega_0}{\omega_p})\hat\kappa\cdot\hat\kappa$ and $\mu_{eff}(\frac{\omega_0}{\omega_p})$ respectively. Figure \ref{example1} (c) (Figure \ref{example2} (c)) compares the  prediction given by the leading order dispersion relations $dk=\sqrt{n_{eff}^{2}}\frac{\omega_0}{\omega_p}$ (solid lines) with  the numerical approximation of the dispersion relation given by black dots.  Although the length scale of the microstructure is not infinitesimally small compared to the radiation wavelength, the numerically calculated points (black dots) fall near the solid lines predicted by the leading order dispersion relation. Notice that the magenta area is the prediction of the band of double negative leading order effective properties while the green area is for the double positive band.  Both graphs show that the leading order dispersion relation given by the power series can well predict the Bloch  wave modes in the actual crystal up to about $20\%$ smaller than the wavelength.

The simulations show that the frequency range spanned by the double negative propagation bands increases as the thickness of the plasmonic coating decreases. 

\section{Acknowledgments}
This research is supported by NSF grant DMS-0807265 and AFOSR grant FA9550-05-0008.

\end{document}